\documentclass{amsart}

\usepackage{amsmath, amsthm, amssymb}
\usepackage{graphicx}
\usepackage{multicol}
\usepackage{footmisc}
\usepackage{newtxtext}
\usepackage{newtxmath}

\title{Binary quadratic forms: modern developments}
\author{Ayberk Zeytin}
\address{Department of Mathematics, Galatasaray University}
\email{azeytin@gsu.edu.tr}

\usepackage{caption, subcaption, amsfonts}
\usepackage{epstopdf}
\usepackage{relsize}
\usepackage{pgf, float}
\usepackage{hyperref}

\usepackage{tikz}
\usetikzlibrary{matrix,arrows}
\usetikzlibrary{decorations.pathmorphing,mindmap}
\usetikzlibrary{mindmap,matrix,arrows,decorations.pathmorphing,patterns,fadings}

\newcommand{\F}{\mathcal{F}}
\newcommand{\QQ}{\mathbf{Q}}
\newcommand{\ZZ}{\mathbf{Z}}
\newcommand{\CC}{\mathbf{C}}
\newcommand{\RR}{\mathbb R}
\newcommand{\HH}{{\mathcal{H}}}
\newcommand{\aut}{\mathrm{Aut}}
\newcommand{\stab}{\mathrm{Stab}}
\newcommand{\ITEM}{$\blacktriangleright$}
\newcommand{\G}{\mathcal{G}}
\newcommand{\pslz}{\mathrm{PSL}_2 (\ZZ)}

\newcommand{\ol}[1]{\overline{#1}}

\newcommand{\sm}{\setminus}

\newtheorem{theorem}{Theorem}[section] 

\newtheorem{proposition}[theorem]{Proposition}

\newtheorem*{definition}{Definition}

\makeindex             

\begin{document}
	

%
%
\maketitle

\begin{abstract}
	In this chapter, we offer a historical stroll through the vast topic of binary quadratic forms. We begin with a quick review of their history and then an overview of contemporary algebraic developments on the subject.
\end{abstract}

\section{Introduction}

It seems that the first encounter of mankind with binary quadratic forms\footnote{A binary quadratic form is a homogeneous element of the ring $\ZZ[x,y]$ of degree two.} is through the famous Pythagorean triples, as we will see below. In fact, the Babylonian tablet Plimpton 322, dating back to around 1800 BCE, represents the earliest documented interaction between humans and binary quadratic forms. As far as the author is aware, there is a gap until Euclid, who lived around 300 BCE. Euclid studied quadratic irrationalities, which may be viewed as an indirect study of the topic. Diophantus of Alexandria (circa 3rd century CE) made some early investigations into quadratic forms in his work "Arithmetica". Namely, Diophantus dealt with special cases of quadratic equations and their solutions. The Indian mathematician Brahmagupta (circa 598 CE - 668 CE) made significant contributions to the study of indeterminate equations, including quadratic ones. For instance, he used the famous formula for the roots of quadratic equations. 

Then, until Fermat (1607 - 1665), no significant contributions to the topic exists to our knowledge. He worked not only on Pythagorean triples but also on some equations of Pell type. It seems that binary quadratic forms were studied heavily during and after the 18th century, especially with contributions from many prominent mathematicians, including Euler (1707 - 1783), Lagrange (1736 - 1813), Legendre (1752-1833), Gauss (1777 - 1855), and Dirichlet (1805 - 1859). One addition that we propose to this list is Poincaré\footnote{The mainstream LLMs do not recognize Poincaré in a list of mathematicians who made contributions to the topic of binary quadratic forms. In fact, here is an output of ChatGPT, V3.5 when asked about Poincaré's contributions: "Henri Poincaré is not primarily known for his work on binary quadratic forms. In fact, the study of binary quadratic forms is a field of mathematics that was extensively developed by other mathematicians, such as Gauss, Lagrange, and Jacobi."}\footnote{This text was finalized before the appearance of \cite{goldstein/poincare/and/arithmetic}, to which the reader is encouraged to consult for a dedicated and complete treatise.}, whose contributions are mostly overlooked.

The theory of binary quadratic forms can be viewed, at least historically, as the birth of algebraic number theory. Indeed, it was the observation that $X^{2} + Y^{2}$ is enough for representing primes that are congruent to 1 modulo 4, whereas a similar statement has to include not only $X^{2} + 5Y^{2}$ but also $2X^{2} + 2XY +3Y^{2}$, that gave rise to many questions. Such problems led to the definition of classes of binary quadratic forms. Gauss carried this idea further and defined a binary operation, called Gauss product, on the set of classes of primitive binary quadratic forms of discriminant $\Delta$, denoted by $\F(\Delta)$. It turned out that  $\F(\Delta)$ endowed with the Gauss product is an abelian group which is isomorphic to the narrow class group of a corresponding quadratic number field (whenever $\Delta$ is a fundamental discriminant) in our modern terminology. Contemporary algebraic number theory is so well developed that its one of the many raison d'\^{e}tre, namely the representation problem of binary quadratic forms, is not mentioned quite often and is not visible to an inexperienced eye. In this work, we focus more on the modern algebraic developments of the binary quadratic forms, leaving many topics untouched, for instance, the analytic aspects and, in particular, the theta functions of definite and indefinite binary quadratic forms. Reader interested in such topics is invited to consult \cite{vigneras/series/de/theta/des/formes/quadratique/indefinies} and \cite{zwegers/thesis} as a starting point. 

The organization of the paper is as follows: The next section is a quick review of the history. This part is neither complete nor detailed. We take refuge in the existence of many wonderful historical accounts on number theory (e.g.\cite{weil/number/theory, dickson/history/1,dickson/history/2,dickson/history/3}). However, we give a little bit detailed account about Poincar\'{e}'s contributions to the field, which seem to have been overlooked. The last section is devoted to some central and contemporary developments concerning the algebraic aspects of indefinite binary quadratic forms. 

\section{A very brief history of binary quadratic forms}

Babylonian tablet Plimpton 322, dated nearly as 1800 BCE, is the first known documented encounter between mankind and binary quadratic forms\footnote{The tablet was found in Iraq and is now at Columbia University, see \texttt{http://www.math.ubc.ca/{$\sim$}cass/courses/m446-03/pl322/pl322.html}.}. It contains the \emph{first} list of 15 Pythagorean triples, that is integral solutions of the equation $X^{2}+Y^{2}=Z^{2}$ with a couple mistakes. The list contains the triple $(12709, 13500, 18541)$ which appears to be the largest. It turns out that Babylonians knew the \emph{rule} in producing such triplets but interested merely in some particular cases. They have also made certain attempts in the direction of proof,\cite{hoyrup/pythagorean}, yet whether the Babylonians knew the proof or not is still in debate. The book \cite{neugebauer/sachs/mathematical/cuneiform/texts} seems to be the first text where the tablet in question is first analyzed. Reader interested in Babylonian mathematics is invited to consult \cite{neugebauer/zur/geometrischen/algebra} and references therein. 

Book X of Euclid's \emph{Elements} (circa 300 BCE) is reserved for the study of quadratic irrationalities. This is the book where Euclid describes his famous algorithm (Proposition~3). Though it is closely related, Euclid's viewpoint did not pay much attention to binary quadratic forms. We refer to \cite{knorr/elements} for a treatment of Book X of Elements in modern terms.

The Cattle problem, attributed to Archimedes (approximately around 250 BCE), is another instance where we encounter binary quadratic forms. For the exact statement of the problem, the interpretation of the problem given by Wurm (sometimes called \emph{Wurm's problem}) and its reduction to the equation:
$$1 + 4729494 u^{2} = t^{2}; $$
(where $2|u$ and $4657|u$) can be found in \cite{vardi/cattle/problem}.  There was a small error in the solution of Amthor, \cite{amthor/cattle/problem}, which came as late as 1880. This mistake was then corrected in \cite{grosjean/cattle/problem} more than 100 years later! 

Some place for binary quadratic forms has also been reserved in ``Arithmetica'' of Diophanthus of Alexandria\footnote{One should mention that the exact dates of birth and death of Diophantus are still in dispute, although we know, from an old riddle composed by Metrodorus, that he died when he was 84. The exact riddle reads: ``This tomb holds Diophantus. Ah, how great a marvel! the tomb tells scientifically the measure of his life. God granted him to be a boy for the sixth part of his life; and adding a twelfth part to this he clothed his cheeks with down; He lit him the light of wedlock after a seventh part, and five years after his marriage He granted him a son. Alas! late-born wretched child, after attaining the measure of half his father's life, chill Fate took him. After consoling his grief by this science of numbers for four years, he ended his life.''\cite[pg. 93]{riddle/diophantus}. It is also arguable, though probably/apparently unnecessary, that the name Diophantus is in fact an ancient analogue of Bourbaki, i.e. a collective author, \cite{schappacher/diophant}.}(3rd century CE). This is a book comprised of 130 problems. There, once again, we meet Pythagorean triples, with a somewhat more systematic treatment. In fact, starting from two Pythagorean triples, Diophantus knew how to derive a third one using the identity
\begin{eqnarray}
(X^{2} + Y^{2})(X'^{2} + Y'^{2}) = (X X' \pm YY')^{2} + (X Y' \mp X'Y )^{2}.
\label{eq:pythagorean}
\end{eqnarray}

This is related, in particular, to problems 8\footnote{To split a given square in two squares.} and 9\footnote{To split a given number which is the sum of two squares into two other squares.} of Book II of Arithmetica. Indeed, Diophanthus refers to his book ``Porismata'' (auxiliary theorems) quite frequently in ``Arithmetica'' and according to Weil 'it is perhaps not too far-fetched to imagine that (\ref{eq:pythagorean}) may have been one of his porismata',\cite[pg.~11]{weil/number/theory}. More on the history of sums of two squares can be found in \cite[Chapter~6]{dickson/history/2}, as well.

As far as the author is aware, there is a gap in the development of binary quadratic forms until the time of Fermat except for Leonardo of Pisa (1170 - 1250), who is rather known by the name Fibonacci. His \emph{Liber Quadratorum} from 1225 (see \cite{sigler/fibonacci} for an English translation) concerns mostly with square numbers and miscellaneous problems around squares, starting with the famous problem of the representability of integers by sums of two squares. It was Euler who proved relevant statements for the forms $x^{2} + 2y^{2}$ and $x^{2} + 3y^{2}$. For a trained eye, the generalization of (\ref{eq:pythagorean}) to 
\begin{eqnarray*}
(X^{2} + NY^{2})(X'^{2} + NY'^{2}) = (X X' \pm NYY')^{2} + N(X Y' \mp X'Y )^{2}
\end{eqnarray*}
is quite clear, but we had to wait Euler and his book algebra for this and its variant
\begin{eqnarray*}
(X^{2} - NY^{2})(X'^{2} - NY'^{2}) = (X X' \pm NYY')^{2} - N(X Y' \pm X'Y )^{2}
\end{eqnarray*}
to appear. Interested reader is invited to consult \cite{konen/geschichte} and \cite{whitford/pell/equation} for further details on the history of the Pell equation. Pell's equation has other applications in solving other Diophantine equations, see for instance, \cite[Chapter 8]{mordell/diophantine/equations}. 

Similar ideas did not lead to a solution for the case $x^{2} + 5y^{2}$. This was known to Fermat. Euler had also failed to explain this phenomenon. It was Lagrange who came up with the idea to separate the $\pslz$ classes of binary quadratic forms of discriminant $D$ according to the classes they represent in $(\ZZ/|D|\ZZ)^{\times}$. Legendre(1752-1833) took the studies of Lagrange further and generalized what Lagrange achieved. Indeed, the identity 
$$ (2x^{2} + 2xy + 3y^{2}) (2x'^{2} + 2x'y' + 3y'^{2}) = X^{2} + 5Y^{2}; $$
where $X = 2xx' + xy' + yx' + 3yy'$ and $Y = xy'-yx'$ is used by Lagrange to prove quadratic reciprocity concerning the binary quadratic form $X^{2} + 5Y^{2}$. This identity is in fact one of the first examples of composition of two binary quadratic forms. Legendre recognized that similar identities can be reproduced for any two binary quadratic forms having the same discriminant laying the foundations of composition. The missing pieces of this composition that prevented the set of binary quadratic forms of fixed discriminant being a group (including the composition as defined by Legendre being multi-valued) has been completed by Gauss in his \emph{Disquisitiones Arithmeticae}, \cite{disquisitiones}. This corner stone contribution - the author believes - is the first work where binary quadratic forms are treated not as tools to solve other problems but rather as objects of interest. Gauss had also treated many other concepts such as \emph{ambiguous forms}\footnote{A binary quadratic form $f = (a,b,c)$ is called ambiguous if there is some $\gamma \in \pslz$ so that $\gamma\cdot f = (a,-b,c)$.} and reduced forms. Gauss had computed class groups of certain imaginary quadratic number fields and put forward his famous class number problems. One should also mention Dirichlet, especially \cite{dirichlet/sur/usage/des/series/infinies,dirichlet/recherches/sur/diverses/applications}, where the famous class number formula for imaginary quadratic number fields is given. 

To the above list of prominent mathematicians, we would like to propose an addition : Poincar\'{e} (1854 - 1912). Poincar\'{e}'s works are not only on binary but also on ternary and other quadratic as well as cubic forms of several variables. In what follows, we offer a look at his contributions to only binary quadratic forms in chronological order.  

The first work of Poincar\'{e} concerning quadratic forms is ”Sur quelques propriét\`{e}s des formes quadratiques”, \cite{poincare/1879a} where he announces that to each form he can associate a corresponding complex number defined as a series. Four properties of these complex numbers are listed. The first one claims that they can be computed by evaluating definite integrals. The remaining three properties are announced so as to attack the classical problems pertaining to quadratic forms; namely, the solution of the corresponding Pell equation, equivalence problem : determining whether two forms are equivalent, and minimum problem : determining the smallest positive integer represented by the quadratic form. These problems were already addressed by Gauss in \cite{disquisitiones}. More precisely, in \S~198 the Pell equation\footnote{Given a form $(M, N,P)$ whose determinant is $D$, and $m$ the greatest common divisor of the numbers $M$, $2N$,  $P$: to find the smallest numbers $t$, $u$ satisfying the indeterminate equation $t^{2}-Du^{2} = m^{2}$} and its solution is discussed, in \S~158 the equivalence problem\footnote{ "Given any two forms having the same determinant, we want to know whether or not they are equivalent, whether or not they are properly or improperly equivalent, or both (for this can happen too).} is discussed and in \S~172 the definition of reduced forms is given which is closely related to the minimum problem. 

In the continuation of this work, \cite{poincare/1879b}, his work on the equivalence problem is announced. Here, Poincar\'{e} mentions the lattice corresponding to a given definite quadratic form. Given two equivalent forms, he mentions that the associated lattices differ by a rotation of angle $\theta$, which he knows how to compute. He then introduces the series 
\begin{align}
\sum_{(x,y) \in \ZZ \times \ZZ \sm \{(0,0)\}} \frac{1}{(ax^{2} + bxy + cy^{2})^{k}}
\label{eq:invariants/arithmetiques}
\end{align}
as an arithmetic invariant\footnote{According to Poincar\'{e} an arithmetic invariant of a form is a quantity invariant under the action of $\mathrm{SL}_{2}(\ZZ)$.} of a given definite binary quadratic form $aX^{2} + bXY +cY^{2}$. The fact that such an invariant can only be a necessary condition is immediately noted.  This gap is filled by defining what is called a covariant function and two such examples are given. It must be noted that the series given in (\ref{eq:invariants/arithmetiques}) appear also in the above mentioned works of Dirichlet, e.g. \cite[\S~7]{dirichlet/sur/usage/des/series/infinies}, as well.  In addition, Poincar\'{e} announces that using these methods he can write a prime number which is congruent to 1 modulo 4 as a sum of two squares. 

The article \cite{poincare/1880} is where Poincar\'{e} uses his correspondence between binary quadratic forms (both definite and indefinite) and lattices and interprets known facts in terms of the geometry (mostly fundamental triangles) of the lattices. The first two sections are devoted to developing algebraic operations (addition, multiplication) on lattices in $\CC$. In the following section, given a binary quadratic form $f = aX^{2} + bXY +cY^{2}$, the lattice generated by $ (b/\sqrt{a}, \sqrt{|b^{2}-4ac| / Da}) $ and $ (\sqrt{a},0) $ is associated to the form; where $D$ is the square-free part of the discriminant of $f$. Various facts concerning this correspondence is discussed. In addition, the reduction algorithm (both for definite and indefinite forms) is also interpreted - using fundamental triangles - in terms of the theory developed in first two sections.  The next section is devoted to defining a group operation which is a geometric version of the group operation defined by Gauss, see the note following this article in \cite{oeuvres/poincare/5} by A.~Ch\^{a}telet. The last section of this work is related to Kummer's (and Dedekind's) theory of ideal numbers. Here, Poincar\'{e} explains how one deals with factorization of ideals in prime ideals in terms of lattices. 

In his work "Sur les invariants arithmetiques"\footnote{Poincar\'{e} published another article with the same title in 1905, \cite{poincare/sur/les/invariants/1905} which will be discussed this below.}, \cite{poincare/sur/les/invariants/1882}, Poincar\'{e} describes a series associated to a linear form. Then, he explicitly describes the series (\ref{eq:invariants/arithmetiques}) that he announced in \cite{poincare/1879b} in addition to another series denoted by $\theta$ by writing the form as a product of two (conjugate) linear forms.  Poincar\'{e} argues (without proof) that these series admit an integral representation, as well. Poincar\'{e} ends this paper with the remark that both for definite and indefinite binary quadratic forms, the series he introduced can be used to distinguish the equivalence of two binary forms.

Poincar\'{e}'s works on quadratic forms focuses on generalization of various known results to forms of more variables and/or of higher degree, for instance, \cite{poincare/sur/les/formes/cubiques/tern/quatern,poincare/sur/la/reproduction,poincare/fonctions/fuchsiennes/et/arithmetique} between 1882 and 1905. This does not mean that Poincar\'{e} did not think about such problems, but rather that he thinks of problems related to arithmetic from a mostly function theoretical viewpoint\footnote{This was remarked by Poincar\'{e}, as well.}. His 1905 article entitled "Sur les invariants arithmetiques", \cite{poincare/sur/les/invariants/1905}, is a very good example of this viewpoint. Indeed, after discussing arithmetic invariants of linear forms, Poincar\'{e} uses this theory to develop similar invariants for quadratic forms. The definite and the indefinite cases are treated separately. Poincar\'{e} interprets and also extends Dirichlet's results using  the theory developed.  We refer to \cite{bergeron/les/invariants/arithmetiques/de/poincare} for a detailed account of this article using modern terms. 

\section{Modern Era}

In this section, we will concentrate on the more modern developments - by modern, we mean the main contribution in the last 50 years. More precisely, we will talk about the work of Bhargava where the Gauss product of two forms is interpreted in terms of what is called Bhargava cubes inspired from the group law on elliptic curves, Zagier's alternative reduction theory, Penner's geometric interpretation of Gauss product for definite binary quadratic forms  and the contribution of the author where the class of an indefinite binary quadratic form is interpreted as a bipartite ribbon graph with specific features. 

\subsection{Notation and terminology}  For an element $f(x,y) = ax^{2} + bxy + cy^{2} \in \ZZ[x,y]$ we will write $f = (a,b,c)$, for short. The coefficients $a$, $b$ and $c$ will be viewed - by abuse of notation - as functions on forms, that is, for $f= (a,b,c)$ we have  $a(f) = a$, $b(f) = b$ and $c(f) = c$. We define the \emph{discriminant} of $f$ to be $\Delta(f) := b^{2} - 4ac$, which is always assumed to be non-zero and non-square, hence $f$ cannot be written as a product of two linear forms with integral coefficients. $f$ is called \emph{positive} (resp. \emph{negative}) \emph{definite} if $\Delta(f)<0$ and $a>0$ (resp. $a<0$). If $\Delta(f)>0$; $f$ is called \emph{indefinite}. $f = (a,b,c)$ is called \emph{primitive} if $1$ is the largest integer dividing $a$, $b$ and $c$ at the same time. For any binary quadratic form $f$ and any given integer $N$, the \emph{set of representations} of $N$ by $f$ is defined as :
$$R_{f}(N) := \{ (X,Y) \in \ZZ^{2} \colon f(X,Y) = N \}.$$
Whenever $R_{f}(N) \neq \emptyset$ we say that $N$ is \emph{represented} by $f$. The set of all integers that are represented by $f$ becomes
$$R_{f} = \{N \in \ZZ \colon R_{f}(N) \neq \emptyset\}.$$

The \emph{modular group} $\pslz$ acts on this set by change of variables; namely, for any $f = (a,b,c)$ and $\gamma = \begin{pmatrix}
p & q \\
r & s
\end{pmatrix}
$
we define 
\[ \gamma \cdot f = (f(p,r) , 2apq + b(ps + qr) + 2crs , f(q,s)).\]
The equivalence class of $f = (a,b,c)$ is denoted by $[f] = [a,b,c]$. The above concepts, the discriminant (hence $f$ being positive/negative definite or indefinite), primitivity, the set $R_{f}$ are invariant under this action. In the search for a unique representative of each equivalence class, we define a definite binary quadratic form $f = (a,b,c)$ to be \emph{reduced} if the following inequalities are satisfied :
\[ |b| \leq a \leq c. \]
If further one of the inequalities is an equality then we require $b \geq 0$. We say that an indefinite binary quadratic form $f$ is \emph{G-reduced} if we have :
\[ |\sqrt{|\Delta(f)|} - 2|a|| < b < \sqrt{|\Delta(f)|}. \]
We also say that an indefinite binary quadratic form is \emph{semi-reduced} if $ac <0$. In the case of positive definite forms, it turns out that every equivalence class has a unique representative. However, we will see in Section~\ref{sec:carks} that whenever $f$ is indefinite, the equivalence class $[f]$ has at least two different reduced forms. 

Let us denote the set of forms of discriminant $\Delta$ as $\F(\Delta)$. Gauss has defined a binary operation which will be referred to as \emph{Gauss product}. In modern terms, $\F(\Delta)$ together with Gauss product turns out to be an abelian group. To the experts (see \cite{arndt/zur/theorie/kubischen/formen,cayley/memoire/hyperdeterminants,dedekind/composition/binaren/quadratischer/formen,weber/composition/quadratischer/formen}), the definition of Gauss seemed unmotivated. In Section~\ref{sec:cubes} we will mention Bhargava's new approach to the topic.

\subsection{Minus continued fractions : Zagier reduction}
\label{sec:zagier}
In this section, we discuss Zagier's contribution to the theory. Let us begin with :
\begin{definition}{\cite[p. 122]{zagier/zetafunktionen/quadratische/zahlkorper}}
An indefinite binary quadratic form $f = (a,b,c)$ is called \emph{Z-reduced} if 
\[ a,c>0 \mbox{ and } b>a+c. \]
\end{definition}
This definition is slightly different from the definition given by Gauss. For instance, if $f = (a,b,c)$ is G-reduced, then $ac<0$ is satisfied; whereas for Z-reduced forms, $ac>0$.  Given any indefinite binary quadratic form $f$, the transformation 
\[ S(L^{2}S)^{n} = \begin{pmatrix}
n & 1 \\ -1 & 0
\end{pmatrix} \]
is applied successively; where $n = \lceil \frac{b+ \sqrt{\Delta(f)}}{2a} \rceil$.  This results in an eventually periodic sequence of forms. Recall that given an indefinite binary quadratic form $f = (a,b,c)$, one can recover the periodic part of the standard continued fraction of the root of the equation $ax^{2} + bx + c = 0$ from the period of the form $f$ that appear in the reduction process. Zagier reduced forms recover the mentioned correspondence in the case of minus continued fractions\footnote{The minus continued fraction is the process of expressing a given real number as a sum of the form 
\[
[[a_{0};a_{1},a_{2},a_{3}]] := a_0 - \cfrac{1}{a_1 - \cfrac{1}{a_2 - \cfrac{1}{a_3 - \cfrac{1}{\ddots}}}} ;
\] 
where $a_{1},a_{2},\ldots$ are all elements of $\ZZ_{\geq 2}$.}\footnote{Minus continued fractions enjoy properties both similar to (such as the characterization of a real number $x$ being a rational number being having a finite continued fraction expansion) and different from (such as roots of quadratic equations admitting eventually periodic expansions in the standard case whereas the terms in expansion being equal to 2 eventually in the minus continued fractions case) standard continued fractions.}.

The main application of this theory is an interesting result - known as Hirzebruch theorem - concerning class numbers of imaginary quadratic number fields :

\begin{theorem}{\cite[\S~14,Satz~3]{zagier/zetafunktionen/quadratische/zahlkorper}}
Let $p>3$ be a prime number which is congruent to $-1$ modulo $4$. Write $\sqrt{p} = [[a_{0};\overline{a_{1},a_{2},\ldots,a_{r}}]]$
with the part $ \overline{a_{1},a_{2},\ldots,a_{r}}$ representing the minimal period of the minus continued fraction expansion. If $\QQ(\sqrt{p})$ has narrow class number $1$, then the class number of $\QQ(\sqrt{-p})$ is 
\[ \frac{1}{3}\left(a_{1} + a_{2} + \ldots + a_{r} \right)-r . \]
\end{theorem}
One should consult \cite{smith/minimal/periods/in/Zagier/reduction} for further discussion regarding Z-reduced forms, which completes the analogy concerning the reduction algorithm associated to this viewpoint.

\subsection{Geometry of Gauss product : Penner's work}
\label{sec:geometry/of/composition}
In \cite{penner/geometry/of/gauss/product} Penner discusses a geometric interpretation of the Gauss law. In order to describe such an interpretation one needs to revisit product law on binary quadratic forms as defined by Dirichlet. Let us recall that given two binary quadratic forms $f_{1} = (a_{1},b_{1},c_{1})$ and $f_{2} = (a_{2},b_{2},c_{2})$, we say that $f_{1}$  and $f_{2}$ are \emph{unitable} if $\Delta(f_{1}) = t_{1}^{2}\Delta$ and $\Delta(f_{2}) = t_{2}^{2}\Delta$ - that is when the square-free part of the discriminants agree. In the case when $\Delta(f_{1}) = \Delta(f_{2})$, we further say that $f_{1}$ and $f_{2}$ are \emph{concordant} if $b_{1} = b_{2}$ and if $a_{2}| c_{1}$ and $a_{1}|c_{2}$. That is, if $f_{1} = (a_{1},b,a_{2}c)$ and $f_{2} = (a_{2},b,a_{1}c)$. In such a case, the product $f_{1}\cdot f_{2}$ is defined simply as $(a_{1}a_{2},b,c)$. It is proven by Dirichlet that given two forms, say $Q_{1},Q_{2}$, of discriminant $\Delta$ one can always find forms $f_{1} \in [Q_{1}]$ and $f_{2} \in [Q_{2}]$ so that $f_{1}$ and $f_{2}$ are concordant, see for instance \cite[p.~335]{cassels/rational/quadratic/forms}.  Moreover, one can already see that unitable forms are concodant but the converse is not true.

Given a positive definite binary quadratic form $f = (a,b,c)$, the equation $aX^{2} + bX + c = 0$ has exactly two conjugate roots : $ \omega_{f} = \frac{-b + \sqrt{\Delta(f)}}{2a} \in \HH = \{z \in \CC \,|\, \mathrm{Im}(z)>0\}$ and $ \ol{\omega_{f}} \in \ol{\HH}$. As $a,b,c \in \ZZ$, $\omega_{f} \in \HH \cap \QQ(\sqrt{-1})$. Conversely, given an element $\omega =  \frac{\alpha}{\beta} + \frac{\gamma}{\delta}\sqrt{-1} \in \HH \cap \QQ(\sqrt{-1})$, the integral positive definite binary quadratic form $f_{\omega} = (\beta^{2}\delta, -2 \alpha \beta \delta, \alpha^{2}\delta + \beta^{2}\gamma)$ admits $\omega$ as its corresponding element. When $\gcd(\alpha,\beta)$ and $\gcd(\gamma,\delta)$ are $1$, $f_{\omega}$ is primitive. One computes that $\Delta(f_{\omega}) = -4\beta^{4}\gamma\delta$. 

In the view of the above correspondence, given two different points $\omega_{1}, \omega_{2} \in \HH \cap \QQ(\sqrt{-1})$ the forms $f_{\omega_{1}}$ and $f_{\omega_{2}}$ are unitable if and only if there are constants $\lambda_{1},\lambda_{2} \in \ZZ$ so that ${\Delta(f_{\omega_{1}})}/{\lambda_{1}^{2}} = {\Delta(f_{\omega_{2}})}/{\lambda_{2}^{2}}$. In particular, if $\omega_{1}$ and $\omega_{2}$ are solutions of the equation 
\begin{align}
\frac{\Delta(f_{\omega})}{a(f_{\omega})^{2}} = -\alpha^{2} \quad \mbox{ or } \quad \frac{\Delta(f_{\omega})}{b(f_{\omega})^{2}} = -\beta^{2} \quad \mbox{ or } \quad \frac{\Delta(f_{\omega})}{c(f_{\omega})^{2}} = -\gamma^{2} 
\label{eq:loci}
\end{align}
for some $\alpha, \beta, \gamma \in \ZZ_{+}$, then the forms $f_{\omega_{1}}$ and $f_{\omega_{2}}$ are unitable. The functions mentioned in (\ref{eq:loci}) induce the following equations in $\HH$
\begin{align}
v=\frac{\alpha}{2} \quad \mbox{ or }\quad v = \pm \beta u \quad \mbox{ or } \quad u^{2} + \left(v - \frac{1}{\gamma} \right)^{2} = \frac{1}{\gamma^{2}},
\label{eq:loci/in/H}
\end{align}
respectively. Considered with its hyperbolic structure, these loci become horocycles centered at $\infty \in \partial \HH$, hypercycle (to the imaginary axis) through $0 \in \partial \HH$ and horocycle centered at $0 \in \partial \HH$, respectively. We conclude
\begin{proposition}{\cite[Corollary~8]{penner/geometry/of/gauss/product}}
Given $f_{1} = (a_{1},b_{1},c_{1})$ and $f_{2} = (a_{2},b_{2},c_{2})$, $f_{1}$ and $f_{2}$ are unitable if and only if they lie on one of the loci mentioned in (\ref{eq:loci/in/H}).
\label{prop:unitable/forms/and/hypercycles}
\end{proposition}
We will refer to these loci as hypercycles. Suppose now that two forms $f_{1}$ and $f_{2}$ and concordant, hence unitable. Then by Proposition~\ref{prop:unitable/forms/and/hypercycles} we know that there is a hypercycle, say $h$, containing $\omega_{f_{1}}$ and $\omega_{f_{2}}$. The following result which gives us an insight about the geometry of the product $[f_{1}]\cdot [f_{2}]$:
\begin{theorem}{\cite[Corollary~3.11]{penner/decorated/TM/theory}}
If $f_{1}$ and $f_{2}$ are two positive definite, primitive concordant binary quadratic forms lying on a common hypercycle $h$, then the point $\omega_{[f_{1}]\cdot [f_{2}]} \in h$ is the point closest to the origin.
\end{theorem}

\subsection{Bhargava's cubes: composition laws}
\label{sec:cubes}
In his thesis, \cite{bhargava/thesis}, which is accompanied by four articles, \cite{bhargava/composition/I,bhargava/composition/II, bhargava/composition/III, bhargava/composition/IV}, Bhargava introduced generalization of Gauss product on the set of binary quadratic forms of fixed discriminant. The main tool was elements of the 8 dimensional $\ZZ$-module $\mathcal{C} :=\ZZ^{2} \otimes \ZZ^{2} \otimes \ZZ^{2}$ whose elements are identified with cubes of integers or equivalently with octuples $(a,b,\ldots,h)$ of integers :

\begin{figure}[H]
	\centering
	\begin{tikzpicture}
		\matrix (m) [matrix of math nodes, row sep=2em, column sep=2em] {
			& e &  & f \\
			a &  & b &  \\
			& g &  & h \\
			c &  & d &  \\
		};
		
		\draw [-] (m-1-2) -- (m-1-4);
		\draw [-] (m-2-1) -- (m-2-3);
		\draw [-] (m-3-2) -- (m-3-4);
		\draw [-] (m-4-1) -- (m-4-3);
		\draw [-] (m-2-1) -- (m-4-1);
		\draw [-] (m-1-2) -- (m-3-2);
		\draw [-] (m-2-3) -- (m-4-3);
		\draw [-] (m-1-4) -- (m-3-4);
		\draw[-] (m-2-1) -- (m-1-2);
		\draw[-] (m-2-3) -- (m-1-4);		
		\draw[-] (m-4-1) -- (m-3-2);
		\draw[-] (m-4-3) -- (m-3-4);		
	\end{tikzpicture}
	\label{fig:cube}
	\caption{An octuple of integers arranged as a cube.}
\end{figure}
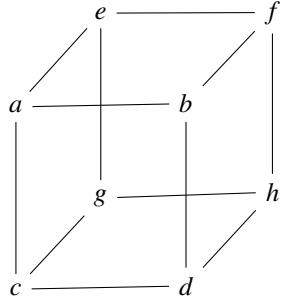

One may partition such a cube into three pairs of $2 \times 2$ matrices using the three linearly independent directions up-down, left-right and front-back :

\begin{align}
\mbox{up-down} &  & 
U = \begin{pmatrix}
	a&e\\b&f
\end{pmatrix} & & D = \begin{pmatrix}
	c&g\\d&h
\end{pmatrix} \notag \\ 
\mbox{left-right} &  & 
L = \begin{pmatrix}
	a&c\\e&g
\end{pmatrix} & & R = \begin{pmatrix}
	b&d\\f&h
\end{pmatrix}
\label{eq:from/cubes/to/forms}
\\ 
\mbox{front-back} &  & 
F = \begin{pmatrix}
	a&b\\c&d
\end{pmatrix} & & B = \begin{pmatrix}
	e&f\\g&h
\end{pmatrix} \notag
\end{align}

From such a partition, we may obtain three binary quadratic forms by declaring
\begin{align*}
f_{UD} = -\det(Ux + Dy) \quad, \quad 	f_{LR} = -\det(Lx + Ry) \quad , \quad 	f_{FB} = -\det(Fx + By) 
\end{align*}

For the pairs of matrices $(U,D)$, $(L,R)$ and $(F,B)$ given in (\ref{eq:from/cubes/to/forms}) we have an action of $\pslz$ :
\[  
\begin{pmatrix}
p&q\\r&s
\end{pmatrix} \cdot (M,N) \mapsto (pM + qN,rM+sN)
\]
inducing an action of $(\pslz)^{3}$ on $\mathcal{C}$. Using this action, one can check that $\Delta(f_{UD}) = \Delta(f_{LR}) = \Delta(f_{FB})$, which will be called the discriminant of the cube $C $ giving rise to these forms. We will denote this discriminant by $\Delta(C)$. The discriminant  of a cube is an invariant of the action of $(\pslz)^{3}$. Primitivity\footnote{One may define a cube $C \in \mathcal{C}$ to be primitive if the forms $f_{UD}$, $f_{LR}$ and $f_{FB}$ are primitive. } of a cube is another invariant of this action.  One has the following result which can be verified with a messy but straightforward  computation :
\begin{theorem}
Let $C \in C $ be a cube giving rise to the forms $f_{UD}$, $f_{LR}$ and $f_{FB}$. The product of the three forms $f_{UD}$, $f_{LR}$ and $f_{FB}$ is the identity, that is
\[ f_{UD} \cdot f_{LR} \cdot f_{FB} = f_{I}; \]
where $f_{I}$ is the identity form - $f_{I} = (1,0,-\Delta(C)/4)$ if $\Delta(C) \equiv 0 \pmod 4$, and $(1,1,-(\Delta(C)-1)/4)$ if $\Delta(C) \equiv 1 \pmod 4$.
\end{theorem}

Following the footprints of the geometric description of the group law on elliptic curves, such an identity can be used to define a group operation, which agrees with the Gauss product by its construction. 

The converse is also possible. That is, given any pair $f_{1},f_{2}$ of primitive binary quadratic forms of discriminant $\Delta$ there exists a primitive cube $C \in \mathcal{C}$ of discriminant $\Delta$ with the property that $f_{UD} = f_{1} $ and $f_{LR} = f_{2}$.

\subsection{\c{C}arks: a combinatorial viewpoint} \label{sec:carks}

\begin{figure*}[H]
	\centering
	\includegraphics[scale=1.25]{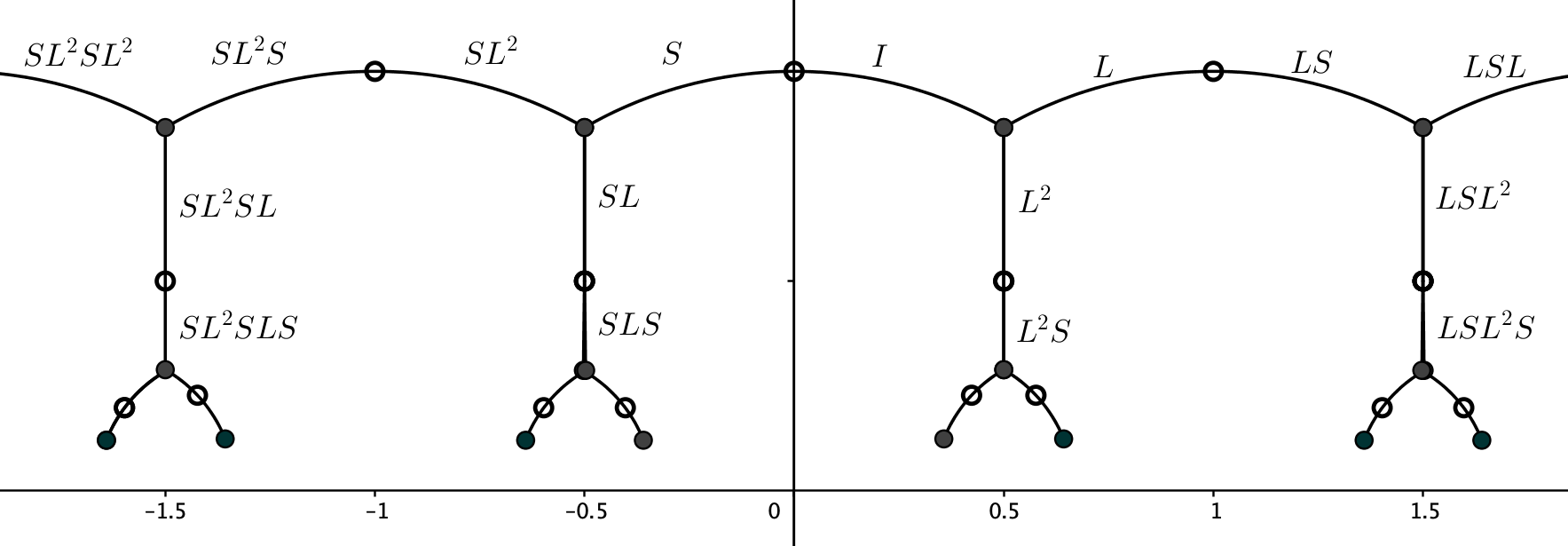}
	\caption{The labeled bipartite Farey tree, $\F$.}
	\label{fig:labeled/bipartite/farey/tree}
\end{figure*}

In recent years, the author and his collaborators, \cite{UZD,reduction}, have given an interpretation of the theory of indefinite binary quadratic forms as a part of the theory of dessins d'enfants. In order to describe the theory, let first fix the elements :
\[ S = \begin{pmatrix}
0&-1\\1&0
\end{pmatrix} \mbox{ and } L = \begin{pmatrix}
1&-1\\1&0
\end{pmatrix} \]
which generate the group $\pslz$ freely. Consider the Farey tree, denoted by $\F$, whose topological realization can be described as follows : let $\bullet$ denote the fixed point of $L$ (that is $\exp(2 \pi \sqrt{-1}/6)$), $\circ$ denote the fixed point of $ S $ (that is, $\sqrt{-1}$) and let $g$ denote the geodesic in $\HH$ joining $\bullet$ to $\circ$. The set of \emph{black} (resp. white) vertices of $\F$ is the orbit of $\bullet$ (resp. $\circ$) and the set of edges is the orbit of $g$ under $\pslz$, see Figure~\ref{fig:labeled/bipartite/farey/tree}. The result of such a construction is a bipartite plane tree where each edge $W\cdot g$ is labeled  with the element $W$. 

The tree $\F$ is endowed with the natural translation action of $\pslz$ defined as left multiplication, which is also compatible with the usual action of $\pslz$ on $\HH$. The automorphism group of any binary quadratic form $f = (a,b,c) $, denoted by $\aut(f)$, is defined as those elements of $\pslz$ which act trivially on $f$, i.e $\stab_{\pslz}(f)$. 

\begin{theorem}{\cite[Theorem~6.12.4]{bqf/vollmer}}
Given any indefinite binary quadratic form $f = (a,b,c)$ the group $\aut(f)$ is generated by one hyperbolic element in $\pslz$, that is $\aut(f) \cong \langle W_{f} \rangle$.
\end{theorem}

\begin{definition}
Given an indefinite binary quadratic form $f$, the infinite bipartite planar graph $\F / \aut(f) = \F/\langle W_{f} \rangle$ is defined as the \c{c}ark\footnote{The word \c{c}ark is Turkish and means wheel. It is pronounced as chark.} associated to $f$ and is denoted by $\G_{f}$.
\end{definition}

\begin{figure}[H]
\noindent \includegraphics[scale=0.75]{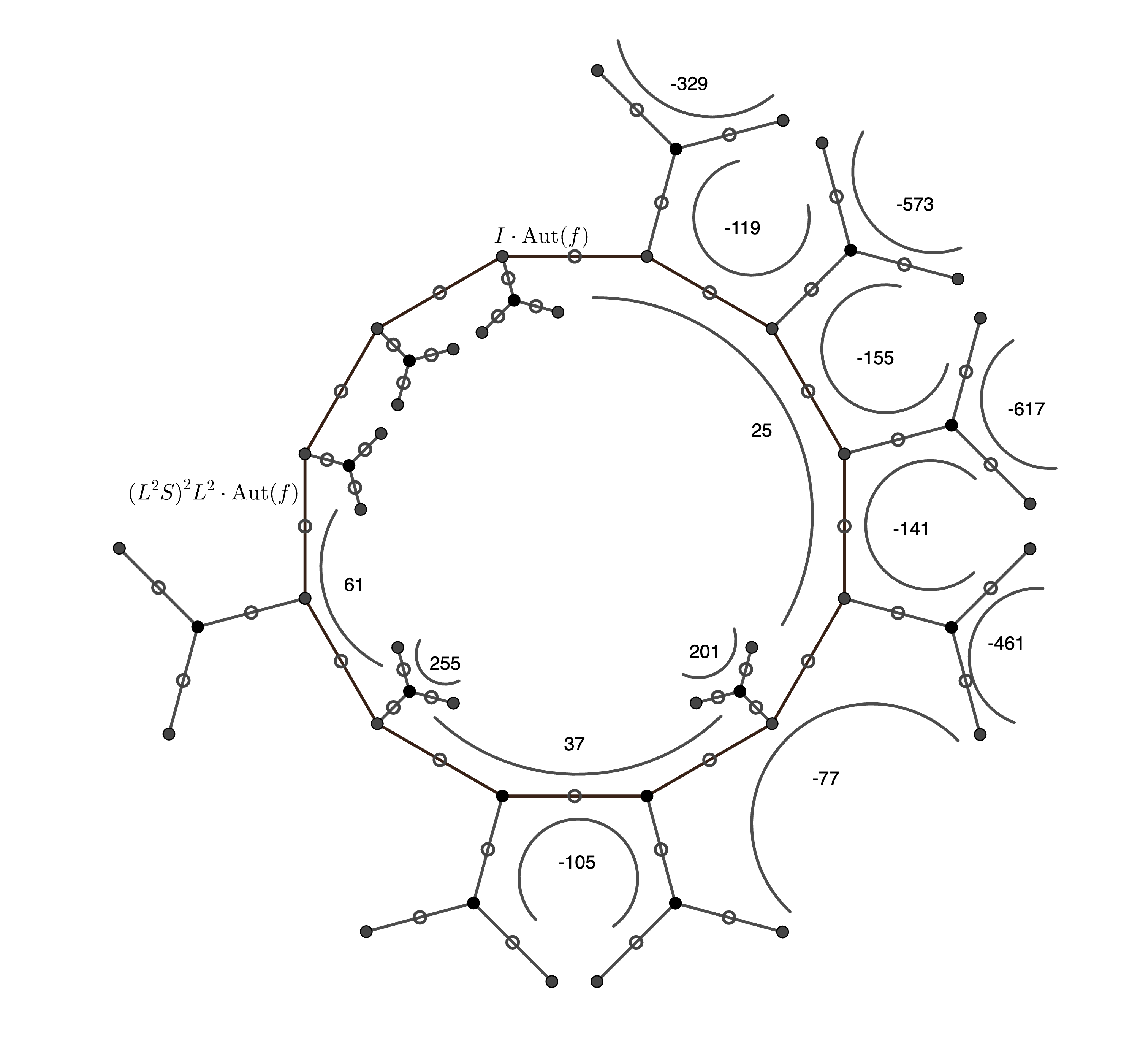}
\caption{The \c{c}ark corresponding to $f = (25,111,-33)$ which is of discriminant $15621 = 3\cdot 41 \cdot 127$.}
\label{fig:cark/example}
\end{figure}

\begin{figure}[H]
	\centering
	\includegraphics[scale=0.40]{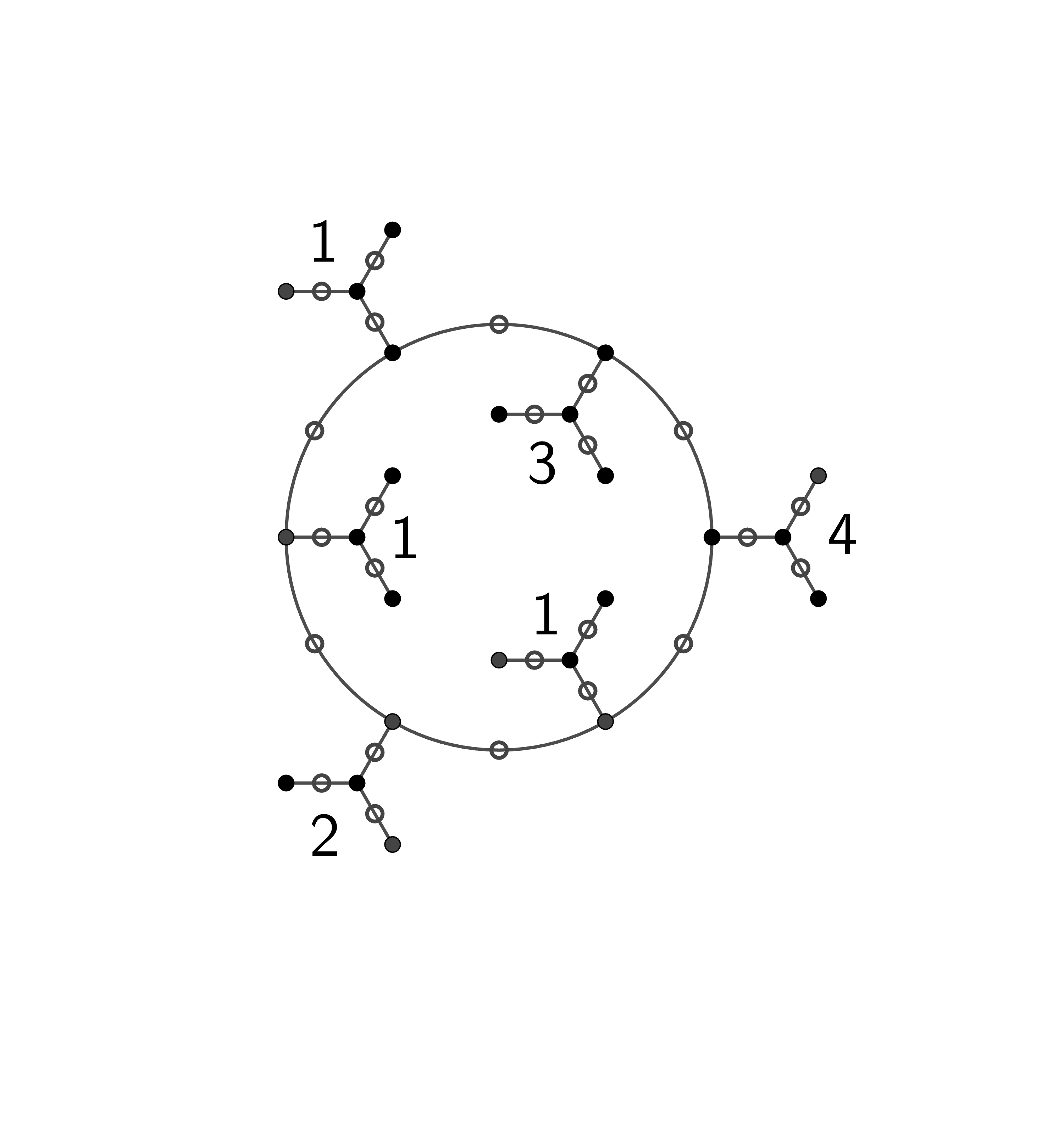}
	\quad
	\includegraphics[scale=0.5]{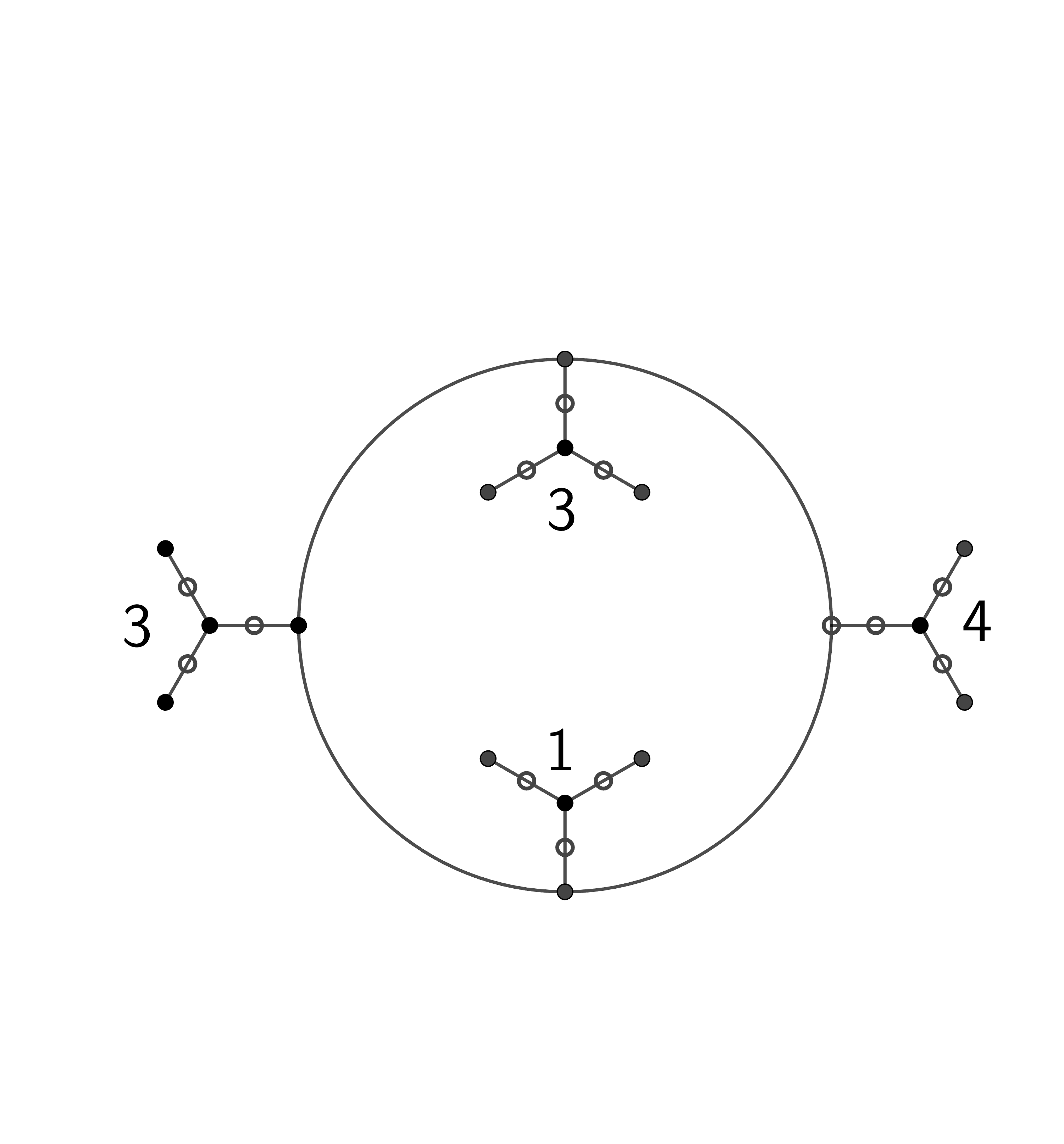}
	\caption{The \c{c}ark corresponding to $f = (25,111,-33)$ and $\mathbf{J}(f)$, respectively.}
	\label{fig:farey/bunches/and/jimm}
\end{figure}

As a result of its construction, the \c{c}ark $\G_{f}$ associated to $f$ enjoys many important features, all of which are explained in Figure~\ref{fig:cark/example} for the form $f = (25,111,-33)$ whose automorphism group is generated by 
\[ W_{f} = (L^{2}S)^{3}(LS)(L^{2}S)(LS)^{2}(L^{2}S)(LS)^{4} = \begin{pmatrix}
7&33\\25&118
\end{pmatrix} \]
\begin{itemize}
\item [\ITEM ] each edge of $\G_{f}$ can be labeled with a coset of $\aut(f)$. The edge labeled with $\aut(f)$ may be regarded as a distinguished (or base) edge of $\G_{f}$. This is a one to one correspondence between the conjugacy class of the generator $W_{f}$ of $\aut(f)$ and the set of edges of $\G_{f}$. 
\item [\ITEM ] the conjugacy class of $W_{f}$ admits a one to one correspondence with the equivalence class of $f$. As a result, each edge of $\G_{f}$ can be labeled with a unique element of $[f]$.
\item [\ITEM ] the vertices of type $\bullet$ (resp. $\circ$) on graph $\G_{f}$ are of degree $3$ (resp. $2$). 
\item [\ITEM ] being a bipartite graph, there is a well defined notion of left and right on $\G_{f}$. A path is called a left-turn path if $(n+1)^{\mbox{\small st}}$ edge of the path is to the left of the $n^{\mbox{\small th}}$  A bi-infinite\footnote{A path on a graph is defined as an ordered sequence of adjacent non-backtracking edges. The first edge is called initial and the last edge is called terminal. A path is called bi-infinite if it does not have neither an initial, nor a terminal edge.} path is called a face of the \c{c}ark.
\item [\ITEM ] the graph $\G_{f}$ has a unique cycle\footnote{A path is called closed if initial and terminal edges agree. A cycle on a graph is defined as a finite closed path.} called the \emph{spine} of $\G_{f}$. If a vertex, an edge, or a face is on the spine, then it will be referred to as a spinal vertex, edge or face, respectively. 
\item [\ITEM ] the spine is the boundary between positive (inner) labeled faces and negative (outer) labeled faces. Spine of $\G_{f}$ is the quotient of the river in Conway's topograph by the automorphism of $f$, \cite{conway/sensual/quadratic/form}. This viewpoint is treated in a more geometric way in \cite{hatcher/topology/of/numbers}.
\item [\ITEM ] there are finitely many Farey branches\footnote{A Farey branch is an infinite bipartite rooted tree which is isomorphic as a bipartite ribbon graph to the quotient $\F/\langle L\rangle$.} attached to the spine at each spinal vertex of type $ \bullet$. Consecutive Farey branches that point in the same direction is called a Farey bunch. For instance, the \c{c}ark in Figure~\ref{fig:cark/example} contains 6 Farey bunches containing, 3,1,1,2,1,4 Farey branches, respectively. Notice that one may use the Farey branches in the Farey bunches together with the natural cyclic order of the spine to specify the \c{c}ark. This is merely the periodic part of the continued fraction expansion of the root of $aX^{2} + bX +c = 0$. Farey bunches allows to depict the \c{c}arks in a more compact way, see Figure~\ref{fig:farey/bunches/and/jimm}, called the Farey bunch representation.
\end{itemize}

The theory around indefinite binary quadratic forms can be interpreted in terms of \c{c}arks. For instance, a form is semi-reduced if and only if it is the label of a spinal edge, a form is Gauss-reduced if it is the label of a spinal edge which is adjacent to the root of an inner Farey bunch. As there are at least two such edges in any given \c{c}ark, this implies that there are at least two reduced forms in any equivalence class. The reduction procedure of Gauss is merely the process of moving the base edge to a base edge, \cite{UZD}. The graph $\G_{f}$ is symmetric with respect to its spine if and only if the form $f$ is ambiguous.

The representation problem admits an algorithmic solution, \cite{reduction,infomod/paper}, as well. To explain this, let us define distance between two vertices to be the number of edges in between, and the root of a face to be the unique vertex (of type $\bullet$) which is closest to the spine. Then, via an \emph{arithmetic progression} on the labels of faces, one may show that these labels increase in absolute value as the distance between the face and the spine increases. Therefore, given any integer $N$, there are only finitely many possibilities for the faces of $\G_{f}$ to be labeled $N$. Therefore, one may solve the representation problem by a search on the faces of $\G_{f}$.

\paragraph{Actions on \c{c}arks.}
The combinatorial/geometric interpretation provided by  \c{c}arks leads to interesting actions. For instance, given a primitive indefinite binary quadratic form $f$ we let $[a_{1},b_{1},\ldots, a_{k},b_{k}]$ be the encoding of the Farey bunches attached to the spine of $\G_{f}$. Then, the \c{c}ark corresponding to  $f^{-1}$ is encoded by $[b_{k},a_{k},\ldots,b_{1},a_{1}]$. 

Another map that acts on the set of \c{c}arks is the map $\mathbf{J}$ studied as a map from $\RR \cup \{\infty\}$ to itself in \cite{uludag/ayral/jimm}. Indeed, Dyer's outer automorphism of the $\mathrm{PGL}_{2}(\ZZ)$ induces such an involution. Despite being discontinuous at the rationals, this involution satisfies a surprising number of functional equations. Indeed, given any continued fraction $[a_{0};a_{1},a_{2},\ldots]$, with $a_{0} \geq 0$ and $a_{i} \geq 1$, is mapped to the continued fraction $[1_{a_{0}-1},2,1_{a_{1}-2},2,1_{a_{2}-2},\ldots] $; where $1_{n}$ denotes the n-tuple $(1,1,\ldots,1)$. By convention $1_{k}$ is  eliminated if it is in the beginning of the expansion and the expansion $[\ldots,a,1_{-1},b,\ldots]$ is replaced by $[\ldots,a+b-1,\ldots]$. In addition, $1_{0}$ is eliminated from the continued fraction expansion, i.e. $[\ldots,a,1_{0},b,\ldots]$ is $[\ldots,a,b,\ldots]$.  Away from finite continued fractions (i.e. rational numbers) and noble numbers (i.e. those real numbers whose continued fraction expansion eventually becomes $[1,1,1,\ldots]$) the map $\mathbf{J}$ satisfies the functional equation :
\[ \mathbf{J}(-x) = -\frac{1}{x}; \]
using which we may define $\mathbf{J}$ for those real numbers whose continued fraction begins with a negative integer. 

The map $\mathbf{J}$ induces a map on the equivalence classes of binary quadratic forms as it preserves those real numbers having eventually periodic continued fraction expansion. For instance, the form $f = (25,111,-33)$ is mapped onto the form $\mathbf{J}(f) = (-43, 46, 13)$ which is of discriminant $4352 = 2^{8}\cdot 17$. Figure~\ref{fig:farey/bunches/and/jimm} shows the Farey bunch representation of \c{c}arks $f$ and $\mathbf{J}(f)$. 

\begin{figure}[H]
	\centering
	\noindent \includegraphics[scale=0.35]{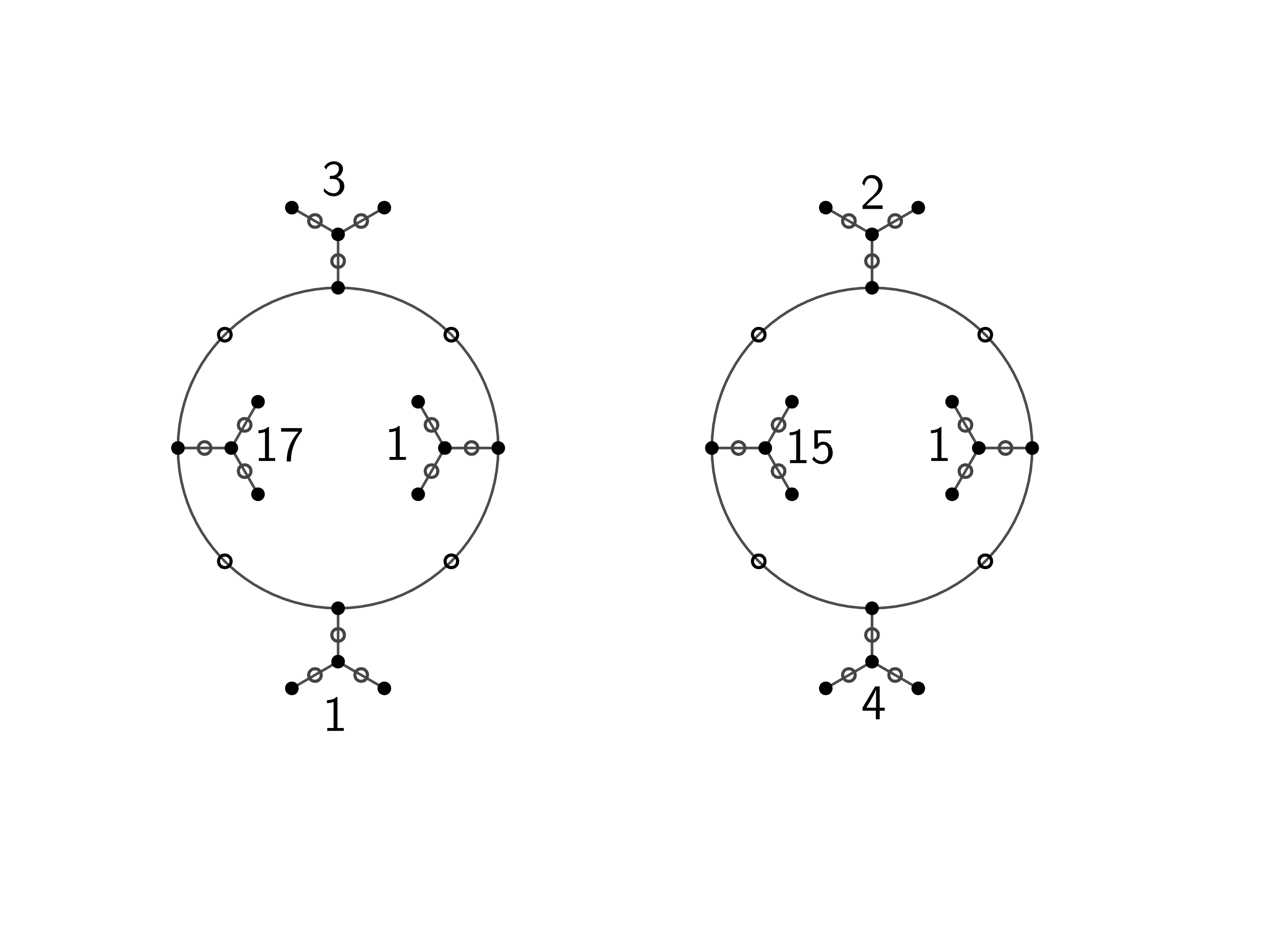}
	\caption{The \c{c}ark corresponding to $g = (-55,89,35)$ and $\mathbf{J}(g)=(-63, 48, 38)$, respectively.}
	\label{fig:cark/15621}
\end{figure}

The form $g = (-55, 89, 35)$ is also of discriminant $15621$ and $g \notin [f]$ as seen in Figure~\ref{fig:cark/15621}. One can quickly compute that $\mathbf{J}(g) = (-63, 48, 38)$ whose discriminant is $11880 = 2^{3}\cdot3^{3}\cdot5\cdot11$

From the above observations, we may see that the action of $\mathbf{J}$ on the class groups can be defined when we consider class groups of all the quadratic orders at once. 

\paragraph{Acknowledgements.} Author is supported by TUBITAK 1001 Research grant 119F405 during the preparation of this manuscript.

\bibliographystyle{alpha}
\def\cydot{\leavevmode\raise.4ex\hbox{.}}

\end{document}